\documentclass{algoritmy}

\usepackage{graphicx}

\usepackage{amsmath}
\usepackage{paralist}
\usepackage{graphics}
\usepackage{graphicx}
\usepackage{epstopdf}
\usepackage{amssymb}
\usepackage{color}
\usepackage{multirow}
\usepackage{booktabs}
\usepackage{bm}
\usepackage{mathrsfs}

\newcommand{\mc}[1]{{\mathcal #1}}
\newcommand{\mb}[1]{{\mathbf #1}}

\title{Evolution of multiple closed knotted curves in space}

\author{Miroslav Kol\'a\v{r}
\thanks{Department of Mathematics, Faculty of Nuclear Sciences and Physical Engineering Czech Technical University in Prague, Trojanova 13, Prague, 12000, Czech Republic}
\and Daniel \v{S}ev\v{c}ovi\v{c}
\thanks{Department of Applied Mathematics and Statistics, Faculty of Mathematics Physics and Informatics, Comenius University, Mlynsk\'a dolina, 842 48, Bratislava, Slovakia. Corresponding author: {\tt sevcovic@fmph.uniba.sk}}
}

\begin{document}

\AlgLogo{1}{10}

\maketitle

\begin{abstract}
We investigate a system of geometric evolution equations describing a curvature and torsion driven motion of a family of 3D curves in the normal and binormal directions. We explore the direct Lagrangian approach for treating the geometric flow of such interacting curves. Using the abstract theory of nonlinear analytic semi-flows, we are able to prove local existence, uniqueness, and continuation of classical H\"older smooth solutions to the governing system of non-linear parabolic equations modelling $n$ evolving curves with mutual nonlocal interactions. We present several computational studies of the flow that combine the normal or binormal velocity and considering nonlocal interaction.
\end{abstract}



\pagestyle{myheadings}
\thispagestyle{plain}
\markboth{M. BENE\v{S}, M. KOL\'A\v{R}, D. \v{S}EV\v{C}OVI\v{C}}{EVOLUTION OF MULTIPLES CLOSED CURVES IN SPACE}

\section{Introduction}
In this work, we focus on the evolution of space curves involving interactions. These one-dimensional structures, which form space curves, are frequently encountered in various scientific and engineering challenges. Connections to dislocation dynamics are discussed in \cite{MBK2010,KKMB2009} along with additional references. Historical research into the dynamics of vortex structures and rings, which align with one-dimensional curves, was initiated by Helmholtz \cite{Helmholtz1858}. The significance of vortex structures in aerospace technology is highlighted in several foundational studies (refer to Thomson \cite{Thomson:67}, Da Rios \cite{Rios:06}, Betchov \cite{Betchov:65}, Arms and Hama \cite{Arms:65}, and Bewley \cite{Bewley:08}). Vortex structures are known to maintain stability over time. This stability is evident in studies of tornadoes and descriptions of volcanic activities (see Fukumoto \emph{et al.} \cite{Fukumoto1987, Fukumoto1991}, Hoz and Vega \cite{Hoz2014}, Vega \cite{Vega2015}). Specific interactions between linear vortex structures, exhibiting dynamic behaviors such as 'frog leaps', are noted (refer to Mariani and Kontis \cite{Mariani2010}). For an overview of vortex dynamics and further discussions on the evolution of closed curves, please see our latest publication \cite{BKS2022} by Bene\v{s}, Kol\'a\v{r}, and \v{S}ev\v{c}ovi\v{c}.

The structure of the paper is as follows. Section 2 revisits the Lagrangian framework for evolving curve families. It introduces a set of evolutionary equations governing the dynamics of interacting curve systems, together with recent findings on the existence and uniqueness of classical H\"older continuous solutions. The proof technique employs the abstract theory of analytic semi-flows in Banach spaces, due to Angenent \cite{Angenent1990, Angenent1990b}. Section 3 concentrates on the numerical discretization approach, utilizing the flowing finite-volume method for spatial derivative discretization and the method of lines to address the resulting ODE systems. Section~6 showcases examples of the dynamics of interacting curves, with interactions shaped by the Biot-Savart nonlocal law, and discusses the development of 3D evolving knotted curves.

\section{Lagrangian description of evolving curves}
We investigate the family of curves $\{\Gamma_t, t\ge 0\}$ evolving in space $\mathbb{R}^3$.  We employ the Lagrangian description of curves, in which the curve is described by a position vector $\mb{X} = \mb{X}(t,u)$ for $t \ge 0$ and $ u \in I $, where $I=\mathbb{R}/\mathbb{Z}\simeq S^1$ is the unit circle. The curve $\Gamma_t$ is then parameterized by
$\Gamma_t = \{ \mb{X}(t,u),  u \in I \}$. The unit tangent vector $\mb{T}$ to $\Gamma_t$ is defined as $\mb{T} = \partial_s \mb{X}$, where $s$ is the unit arc-length parametrization defined by $ds = g du$ where $g=|\partial_u \mb{X}|$ is the relative local length of the curve. Here, $|\cdot|$ denotes the Euclidean norm. The curvature $\kappa$ of a curve $\Gamma_t$ is defined as $\kappa = |\partial_{s}\mb{X}\times\partial^2_{s} \mb{X} | = | \partial_{s}^2 \mb{X} |$. If $\kappa > 0$, we can define the Frenet frame along the curve $\Gamma_t$ with unit normal $\mb{N} = {\kappa}^{-1} \partial_{s}^2 \mb{X}$ and binormal vectors $\mb{B} = \mb{T} \times \mb{N}$, respectively.  
We study a coupled system of evolutionary equations that describes the evolution of closed 3D curves evolving in normal and binormal directions, and the scalar quantity $\varrho$ calculated over the evolving family of curves. More specifically, we focus on the analysis of the motion of a family of $n$ curves evolving in 3D and satisfying the system of geometric equations:
\begin{equation}
\partial_t\mb{X}^i = a^i \partial^2_{s^i} \mb{X}^i + b^i (\partial_{s^i}\mb{X}^i\times\partial^2_{s^i} \mb{X}^i)  + \mb{F}^i +\alpha^i \mb{T}^i,
\qquad \mb{X}^i(\cdot,0)  = \mb{X^i_0}(\cdot), \quad i=1, \ldots, n,
\label{eq:ab}
\end{equation}
which is subject to initial conditions at the origin $t=0$ representing parametrization of the family of initial curves $\Gamma^i_0, i=1, \ldots, n$. Here $a^i=a^i(\mb{X}^i, \mb{T}^i) \ge 0$, and $b^i=b^i(\mb{X}^i, \mb{T}^i)$ are bounded and smooth functions of their arguments, $\mb{T}^i$ is the unit tangent vector to the curve and $s^i$ is the unit arc-length parametrization of the curve $\Gamma^i$. The source forcing term $\mb{F}^i$ is assumed to be a smooth and bounded function. Here
$\mb{F}^i = \mb{F}^i(\mb{X}^i, \mb{T}^i, \gamma^{i1}, \ldots, \gamma^{in})$ is the forcing term and $\gamma^{ij}=\gamma^{ij}(\mb{X}^i, \Gamma^j)$ may depend on the position vector $\mb{X}^i\in\mathbb{R}^3$ and the entire curve $\Gamma^j$. 

As an example of nonlocal source terms $\mb{F}^i, i=1,\ldots,n,$ one can consider a flow of $n=2$ interacting curves evolving in 3D according to the geometric equations:
\begin{equation*}
\begin{split}
\partial_t\mb{X}^1 &= a^1 \partial^2_{s^1} \mb{X}^1 + b^1 (\partial_{s^1}\mb{X}^1\times\partial^2_{s^2} \mb{X}^1 + \gamma^{12}(\mb{X}^1, \Gamma^2),
\\
\partial_t\mb{X}^2 &= a^2 \partial^2_{s^2} \mb{X}^2 + b^2 (\partial_{s^2}\mb{X}^2\times\partial^2_{s^2} \mb{X}^2) + \gamma^{21}(\mb{X}^2, \Gamma^1),
\end{split}
\end{equation*}
where the nonlocal source term
\begin{equation*}
\gamma^{ij}(\mb{X}, \Gamma^j) = \oint_{\Gamma^j}
\frac{(\mb{X}-\mb{X}^j(s^j))\times \mb{T}^j(s^j)}{|\mb{X}-\mb{X}^j(s^j)|^3} ds^j
\end{equation*}
represents the Biot-Savart force measuring the integrated influence of points $\mb{X}^j$ belonging to the second curve $\Gamma^j=\{ \mb{X}^j(u), u\in[0,1]\}$ at a given point $\mb{X}\in\mathbb{R}^3$. 

The tangential velocity $\alpha^i$ that appears in geometric evolution (\ref{eq:ab}) has no impact on the shape of the evolving family of closed curves $\Gamma^i_t, t\ge 0$. This means that the curves $\Gamma^i_t, t\ge 0,$ evolving according to (\ref{eq:ab}) do not depend on a particular choice of the total tangential velocity $v_T^i=\mb{F}^i\cdot\mb{T}^i + \alpha^i$. On the other hand, the tangential velocity has a significant impact on the analysis of evolution of curves from both the analytical and numerical points of view (see e.g., Hou et al. \cite{Hou}, Kimura \cite{Kimura}, Mikula and \v{S}ev\v{c}ovi\v{c} \cite{sevcovic2001evolution, MS2004, MMAS2004}, Yazaki and \v{S}ev\v{c}ovi\v{c} \cite{SevcovicYazaki2012}.  Barrett \emph{et al.} \cite{Barret2010, Barret2012}, Elliott and Fritz \cite{Elliot2017}, investigated the gradient and elastic flows for closed and open curves in $\mathbb{R}^d, d\ge 2$. They constructed a numerical approximation scheme using a suitable tangential redistribution. Kessler \emph{et al.} \cite{Kessler1984} and Strain \cite{Strain1989} illustrated the role of a suitably chosen tangential velocity in numerical simulation of two-dimensional snowflake growth and unstable solidification models. Garcke \emph{et al.} \cite{Garcke2009} applied the uniform tangential redistribution in the theoretical proof of the non-linear stability of stationary solutions for curvature driven flow with triple junction in the plane. 

It is known that the relative local length $g^i/L^i, $ is constant with respect to time $t$, i.e., $\frac{g^i(u,t)}{L(\Gamma^i_t)} =  \frac{g^i(u,0)}{L(\Gamma^i_0)}, \quad u\in I, t\ge 0$ provided that the total tangential velocity $v^i_T$ satisfies $\partial_{s^i} v^i_T= \kappa^i v^i_N  -\frac{1}{L^i} \int_{\Gamma^i} \kappa v^i_N ds^i$ (see, e.g., \cite{Hou}, \cite{Kimura}, \cite{sevcovic2001evolution}, \cite{MS2004}, \cite{MMAS2004}). Another suitable choice of the total tangential velocity $v^i_T$ is the so-called asymptotically uniform tangential velocity proposed and analyzed by Mikula and \v{S}ev\v{c}ovi\v{c} in \cite{MS2004, MMAS2004}. It satisfies $\lim_{t\to \infty} \frac{g^i(u,t)}{L(\Gamma^i_t)} =1$ uniformly with respect to $u\in[0,1]$. This means that the redistribution becomes asymptotically uniform. 

In \cite{BKS2022}, Bene\v{s}, Kol\'a\v{r} and \v{S}ev\v{c}ovi\v{c} generalized methodology and techniques of proofs of the local existence, uniqueness and continuation of solutions from our previous paper \cite{BKS2020} to the case of combined motion of closed space curves evolving in normal and binormal direction and taking into account mutual nonlocal interactions. We proved the result on existence and uniqueness of classical solutions for a system on $n$ evolving curves in $\mathbb{R}^3$ with mutual nonlocal interactions including, in particular, the vortex dynamics evolved in the normal and binormal directions and external force of the Biot-Savart type, or evolution of interacting dislocation loops. 

In the rest of this section we state the result on the existence and uniqueness of classical H\"older smooth solutions to the system of governing equations (\ref{eq:ab}). The method of the proof is based on the abstract theory of analytic semi-flows and the theory of maximal regularity in Banach spaces due to Angenent \cite{Angenent1990, Angenent1990b}. First, we introduce the function space setting. By $h^{k+\varepsilon}(S^1)$ we denote the so-called little H\"older space, i.e. the Banach space which is the closure of $C^\infty$ smooth functions in the norm Banach space of $C^k$ smooth functions defined on the periodic domain $S^1$, and such that the $k$-th derivative is $\varepsilon$-H\"older smooth. Here $0<\varepsilon<1$ and $k$ is a non-negative integer. The norm is given as a sum of the $C^k$ norm and the H\"older semi-norm of the $k$-th derivative. Next, we introduce the scale of Banach spaces of H\"older continuous functions defined in the periodic domain $S^1$:
\begin{equation*}
E_k = (h^{2k +\varepsilon}(S^1) )^3, \quad \mc{E}_{k}= \underbrace{E_{k}\times\ldots\times E_{k}}_{n - times}, \quad k=0,\ 1/2,\ 1
\end{equation*}
We assume that the functions $a^i>0, b^i$, and $\mb{F}^i$ are sufficiently smooth and globally Lipschitz continuous functions (see  \cite[assumtions (H)]{BKS2022} for details). Moreover, we $\alpha^i, i=1, \ldots, n$, is the tangential velocity that preserves the relative local length. 
Assume that the parametrization $\mb{X}_0\equiv (\mb{X}^i_0)_{i=1}^n,$ of initial curves $\Gamma^i_0$ belongs to the H\"older space $\mc{E}_1$,  and it is a uniform parametrization, that is, $|\partial_u\mb{X}^i_0(u)| = L(\Gamma^i_0) >0$ for all $u\in I$ and $i=1,\ldots, n$. 
With regard to \cite[Theorem 4.1]{BKS2022}, there exists $T>0$ and the unique family of curves $\{\Gamma^i_t, t\in[0,T]\}, i=1,\ldots, n$, evolving in 3D according to the system of nonlinear nonlocal geometric equations (\ref{eq:ab}). Their parametrization satisfies $\mb{X}=(\mb{X}^i)_{i=1}^n \in  C([0,T], \mc{E}_1) \cap C^1([0,T], \mc{E}_0)$, and $\mb{X}(\cdot, 0)= \mb{X}_0$. If the maximal time of existence $T_{max}<\infty$ is finite then $\lim_{t\to T_{max}} \max_{i, \Gamma^i_t} |\kappa^i(\cdot, t)| = \infty$.

\section{Flowing finite volumes numerical discretization scheme}

\begin{figure}
\begin{center}
\includegraphics[width=0.3\textwidth]{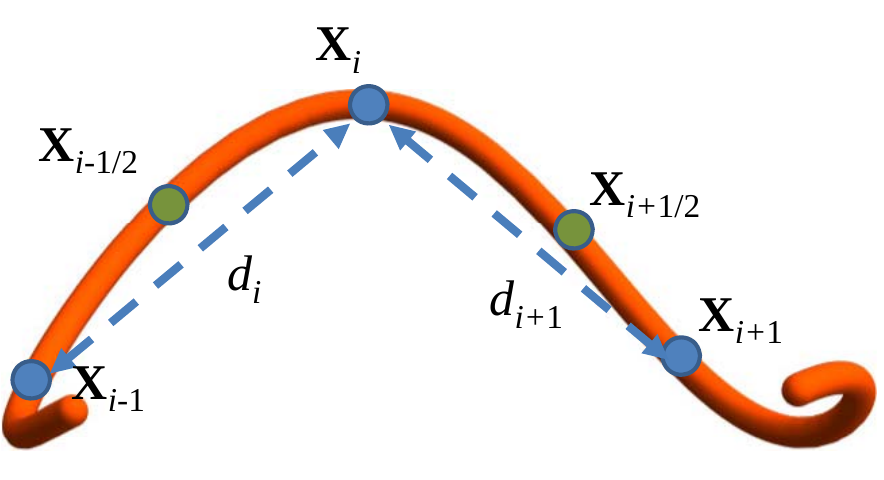}
\end{center}
\caption{Discretization of a segment of a 3D curve by the method of flowing finite volumes.
}
\label{FVMfig}
\end{figure}

In this section, we present a numerical discretization scheme for solving the system of equations (\ref{eq:ab}) enhanced by the tangential velocity $\alpha^i$. Our discretization scheme is based on the method of lines with spatial discretization obtained by means of the finite-volume method. For simplicity, we consider one evolving curve $\Gamma$ (omitting the curve index $i$) and rewrite the abstract form of (\ref{eq:ab}) in terms of the principal parts of its velocity.
\begin{equation}
    \label{eq:num_start}
    \partial_t \mb{X} = a \partial_s^2 \mb{X}
    + b(\partial_s \mb{X} \times \partial_s^2 \mb{X})
    + \mb{F}
    + \alpha \mb{T}.
\end{equation}
We consider $M$ discrete nodes $\mb{X}_k = \mb{X}(u_k)$,  $k = 0,1,2, \ldots, M$ along the curve $\Gamma$. The dual nodes are defined as $\mb{X}_{k \pm \frac12} = \mb{X}(u_k \pm h/2)$ (see Figure \ref{FVMfig}) where $h=1/M$ and $(\mb{X}_k + \mb{X}_{k+1}) / 2$ denote averages in the segments that connect the discrete nodes near them and differs from $\mb{X}_{k \pm \frac12} \in \Gamma$. The $k$-th segment $\mathcal{S}_k$ of $\Gamma$ between the nodes $\mb{X}_{k-1}$ and $\mb{X}_{k}$ represents the finite volume. The integration of equation (\ref{eq:num_start}) in the segment of $\Gamma$ between the nodes $\mb{X}_{k + \frac12}$ and $\mb{X}_{k - \frac12}$ yields
\begin{equation}
\label{eq:num_int}
\begin{split}
\int_{u_{k-\frac12}}^{u_{k+\frac12}} \partial_t \mb{X} |\partial_u \mb{X}| d u = &
\int_{u_{k-\frac12}}^{u_{k+\frac12}} a\frac{\partial}{\partial_u} \left( \frac{\partial_u \mb{X}}{|\partial_u \mb{X}|} \right) d u
+
\int_{u_{k-\frac12}}^{u_{k+\frac12}}
b (\partial_s \mb{X} \times \partial_s^2 \mb{X}) |\partial_u \mb{X}| d u
 \\
& +  \int_{u_{k-\frac12}}^{u_{k+\frac12}} \mb{F} |\partial_u \mb{X}| d u
+ \int_{u_{k-\frac12}}^{u_{k+\frac12}} \alpha \partial_u \mb{X} d u.
\end{split}
\end{equation}
Let us denote $d_k = |\mb{X}_k - \mb{X}_{k-1}|$ for $k=1,2,\ldots,M,M+1$, where $\mb{X}_0 = \mb{X}_M$ and $\mb{X}_1 = \mb{X}_{M+1}$ for closed curve $\Gamma$ and we approximate the integral expressions in (\ref{eq:num_int}) by means of the finite volume method along $\Gamma$ as follows:
\begin{equation*}
\begin{split}
& \int_{u_{k-\frac12}}^{u_{k+\frac12}} \partial_t \mb{X} |\partial_u \mb{X}| d  u
 \approx \frac{d  \mb{X}_k}{d t} \frac{d_{k+1} + d_k}{2},
\\
&\int_{u_{k-\frac12}}^{u_{k+\frac12}}
a \partial_u \left( \frac{\partial_u \mb{X}}{|\partial_u \mb{X}|} \right) d u
\approx
a_k \left(\frac{\mb{X}_{k+1} - \mb{X}_k}{d_{k+1}} - \frac{\mb{X}_k - \mb{X}_{k-1}}{d_k}\right),
\\
&\int_{u_{k-\frac12}}^{u_{k+\frac12}}
b (\partial_s \mb{X} \times \partial_s^2 \mb{X}) |\partial_u \mb{X}| d u
\approx
b_k \frac{d_{k+1} + d_k}{2} \kappa_k (\mb{T}_k \times \mb{N}_k),
\\
&\int_{u_{k-\frac12}}^{u_{k+\frac12}} \mb{F} |\partial_u \mb{X}| d u
\approx  \mb{F}_k \frac{d_{k+1} + d_{k}}{2},
\qquad
\int_{u_{k-\frac12}}^{u_{k+\frac12}} \alpha \partial_u \mb{X} d u\approx  \alpha_k \frac{\mb{X}_{k+1} - \mb{X}_{k-1}}{2}.
\end{split}
\end{equation*}
The approximation of the nonnegative curvature $\kappa$, tangent vector $\mb{T}$ and the normal vector $\mb{N}$, $\kappa \mb{N} = \partial_s \mb{T}$ reads as follows:

\begin{equation*}
\begin{split}
\kappa_k & \approx  
\left|\frac{2}{d_k + d_{k+1}}
\left(\frac{\mb{X}_{k+1} - \mb{X}_k}{d_{k+1}} - \frac{\mb{X}_k - \mb{X}_{k-1}}{d_k}\right)
\right|,
\\
\mb{T}_k & \approx \frac{\mb{X}_{k+1} - \mb{X}_{k-1}}{d_{k+1}+d_k}, \quad
\mb{N}_k \approx \kappa_k^{-1} \frac{2}{d_k + d_{k+1}}
\left(\frac{\mb{X}_{k+1} - \mb{X}_k}{d_{k+1}} - \frac{\mb{X}_k - \mb{X}_{k-1}}{d_k}\right).
\end{split}
\end{equation*}
To discretize the governing equations, we assume that $\partial_t \mb{X}, \partial_u \mb{X}, \mb{F}$ and $\alpha$ are constant over the finite volume between the nodes $\mb{X}_{k + \frac12}$ and $\mb{X}_{k - \frac12}$, taking values $\partial_t \mb{X}_k, \partial_u \mb{X}_k, \mb{F}_k$ and $\alpha_k$, respectively. In approximation $\mb{F}_k$ of the non-local vector-valued function $\mb{F}$, we assume that the curve $\Gamma$ entering the definition of $\mb{F}$ is approximated by the polygonal curve with vertices $(\mb{X}_0, \mb{X}_1, \ldots, \mb{X}_M)$. To find the approximation $\alpha_k$ of the tangential velocity,  we apply a simple trapezoidal integration formula. The values $\alpha_0=\alpha_M$ are chosen so that $\sum_{j=1}^M  \alpha_j d_j =0$. In summary, the semi-discrete scheme for solving (\ref{eq:num_start}) can be written as follows.
\begin{eqnarray}
\frac{d  \mb{X}_k}{d t} \frac{d_{k+1}+d_k}{2}
&=&
a_k \left(\frac{\mb{X}_{k+1} - \mb{X}_k}{d_{k+1}} - \frac{\mb{X}_k - \mb{X}_{k-1}}{d_k}\right)
 + b_k \frac{d_{k+1}+d_k}{2} \kappa_k (\mb{T}_k \times \mb{N}_k)
\label{eq:DCScheme}
\\
&& + \mb{F}_k \frac{d_{k+1} + d_{k}}{2}
 + \alpha_k \frac{\mb{X}_{k+1} - \mb{X}_{k-1}}{2},
 \nonumber
\\
\mb{X}_k(0) &=& \mb{X}_{ini}(u_k), \quad \text{for} \ k = 1, \ldots, M.
\label{eq:DCScheme+1}
\end{eqnarray}
Resulting system (\ref{eq:DCScheme}--\ref{eq:DCScheme+1}) of ODEs  is solved numerically by means of the 4th-order explicit Runge-Kutta-Merson scheme with automatic time stepping control and the tolerance parameter $10^{-3}$ (see \cite{0965-0393-24-3-035003}). We chose the initial time step as $4h^2$, where $h=1/M$ is the spatial mesh size.

\section{Examples of evolution of linked Fourier curves under Biot-Savart external force}
As an example of a non-local source term $\mb{F}$ we consider the external force corresponding to the Biot-Savart law. It represents the integrated influence of all points belonging to the curve $\Gamma=\{ \mb{X}(s), s\in[0,L(\Gamma)]\}$ at a given point $\mb{X}\in\mathbb{R}^3, \mb{X}\not\in \Gamma$. It is given as a line integral:
\begin{equation}
\mb{F}(\mb{X}, \Gamma) = \int_{\Gamma}
\frac{(\mb{X}-\mb{X}(s))\times \partial_s \mb{X}(s)}{|\mb{X}-\mb{X}(s)|^3} ds .
\label{biot-savart-force}
\end{equation}
Let $\Gamma^1$ and $\Gamma^2$ be two non-intersecting closed curves in 3D. The Biot-Savart force is connected with the Gauss linking number and the integral $link(\Gamma^1,\Gamma^2 )$ of $\Gamma^1$ and $\Gamma^2$ can be defined as follows:
\begin{eqnarray*}
&& link(\Gamma^1,\Gamma^2 ) = \frac {1}{4\pi }\oint_{\Gamma^1}\oint_{\Gamma^2}
\frac{det\left(\partial_{s_1} \mb{X}^1(s_1), \partial_{s_2} \mb{X}^2(s_2), \mb{X^1}(s_1) - \mb{X^2}(s_2) \right)}{|\mb{X^1}(s_1) - \mb{X^2}(s_2)|^3} ds_1 ds_2
\\
&&\hskip -0.5truecm = 
- \frac {1}{4\pi }\oint_{\Gamma^1} \mb{F}(\mb{X}^1(s_1), \Gamma^2) \cdot \partial_{s_1} \mb{X}^1(s_1)  ds_1 
=
- \frac {1}{4\pi }\oint_{\Gamma^2} \mb{F}(\mb{X}^2(s_2), \Gamma^1) \cdot \partial_{s_2} \mb{X}^2(s_2)  ds_2 ,
\end{eqnarray*}
where the closed curves $\Gamma^1$ and $\Gamma^2$ are parameterized by $\mb{X}^1(s_1)$ and $\mb{X}^2(s_2)$, respectively. The linking number $link(\Gamma^1,\Gamma^2 )$ belongs to $\mathbb{Z}$. 

\begin{figure}
\begin{center}
\includegraphics[width=0.3\textwidth]{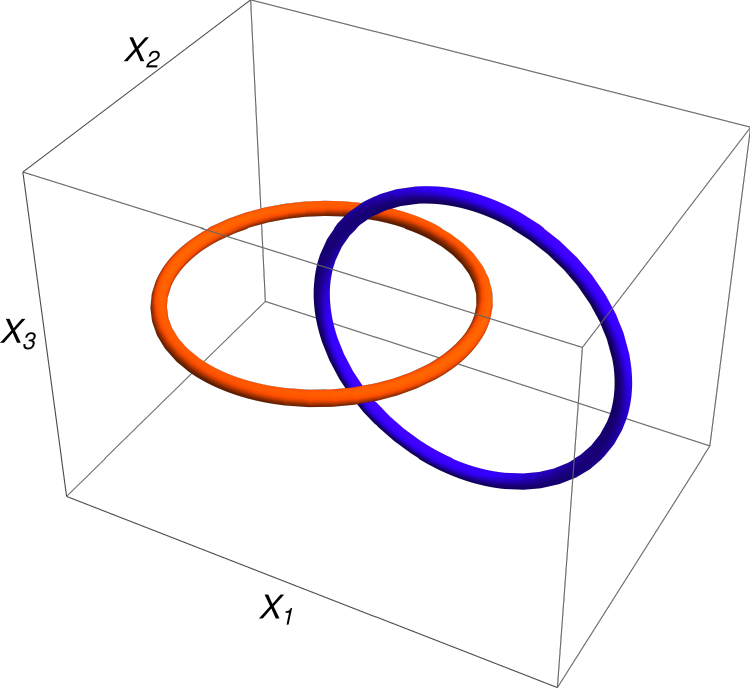}
\qquad
\includegraphics[width=0.25\textwidth]{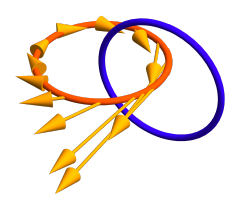}

{\small a) \hskip 5truecm b)}

\caption{Two linked circles a) and the Biot-Savart force vector field induced by $\Gamma^2$ acting on points of $\Gamma^1$, b).}

\label{fig-linked-circles}
\end{center}
\end{figure}

A Fourier curve is a closed curve in 3D that can be parameterized by a finite Fourier series in the parameter $u\in[0,1]$. In Fig.~\ref{fig-linked-circles} a) we show two linked circles $\Gamma^1$ and $\Gamma^2$ with the linking number $link(\Gamma^1,\Gamma^2)=-1$ for parameterization (\ref{eq-crvs1-explicit}), and $link(\Gamma^1,\Gamma^2)=1$ when parameterized by (\ref{eq-crvs2-explicit}):
\begin{eqnarray}
&&\hskip -1truecm\mb{X}^1(u) = (\cos(2 \pi u), \sin(2 \pi u), 0), \ \ \mb{X}^2(u) = (1+\cos(2 \pi u), 0, \sin(2 \pi u)),
\label{eq-crvs1-explicit}
\\
&&\hskip -1truecm\mb{X}^1(u)  = (-\sin(2 \pi u), \cos(2 \pi u), 0),  \mb{X}^2(u)  = (1+\cos(2 \pi u), 0, -\sin(2 \pi u) ).
\label{eq-crvs2-explicit}
\end{eqnarray}
In Fig.~\ref{fig-linked-circles} b) we present two linked circles, and Biot-Savart force vector field induced by $\Gamma^2$ acting on points belonging to $\Gamma^1$. 

\begin{figure}
\begin{center}
\includegraphics[width=0.35\textwidth]{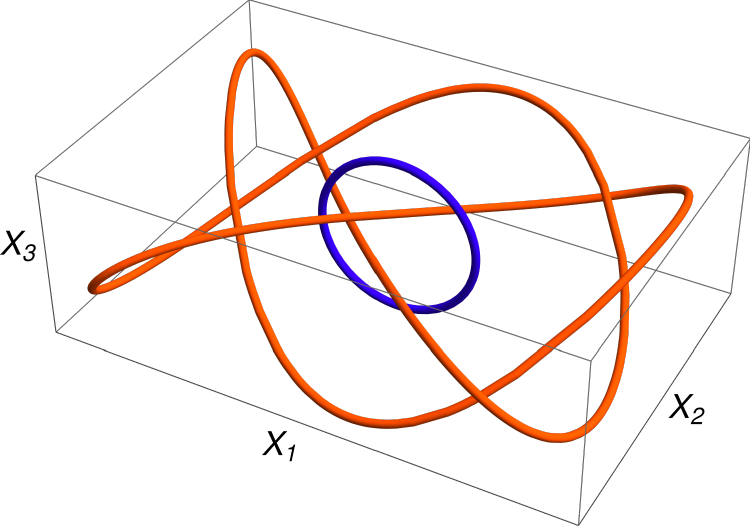}
\qquad 
\includegraphics[width=0.35\textwidth]{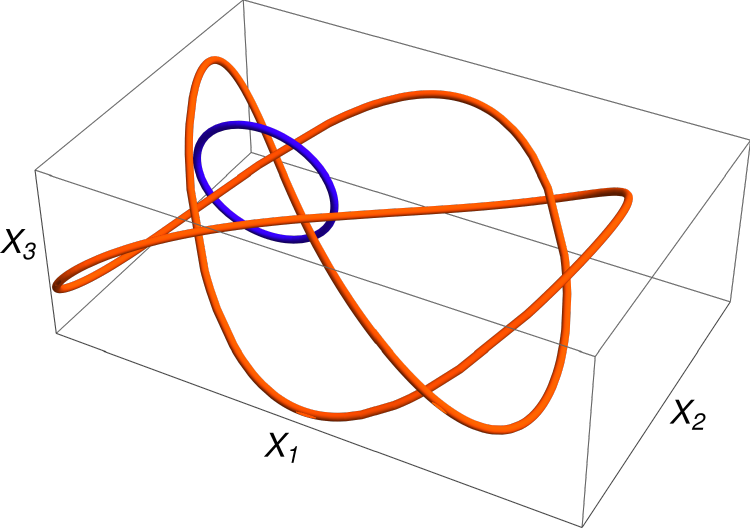}
\centerline{\scriptsize a) \hskip 5truecm b)}
\caption{The Listing's 8-knot curve linked with a circle a) and an ellipse b)}
\label{fig-knot-circles-ellipse}
\end{center}
\end{figure}

The explicit parametrization of the Listing's 8-knot curve $\mb{X}^1$ is given by 
\begin{equation}
\mb{X}^1(u) = (3 \cos(4 \pi u), 2 \sin(6 \pi u + 1/2), (\cos(10 \pi u + 1/2) + \sin(6 \pi u + 1/2))/2 ).
\label{eq-8-curve}
\end{equation}
The Listing's 8-knot curve parameterized by (\ref{eq-8-curve}) is shown in Fig.~\ref{fig-knot-circles-ellipse} a) with a linked-in circle (linking number 0) that is parameterized by: 
\begin{equation*}
\mb{X}^2(u) = (\cos(2 \pi u),\  0\ , \sin(2 \pi u)).
\label{eq-knot-circle1}
\end{equation*}
The Listing's 8-knot curve $\mb{X}^1$  with a linked-in ellipse (linking number -2) that is  parameterized by:
\begin{equation*}
\mb{X}^2(u) = (\cos(2 \pi u)-1.5 , \ 0.5\ , 0.8\sin(2 \pi u))
\end{equation*}
is shown in Fig.~\ref{fig-knot-circles-ellipse} b).

In what follows, we present results of numerical approximation of solutions to the coupled system of governing  PDEs:
\begin{equation*}
\begin{split}
\partial_t\mb{X}^1 &= \partial^2_{s^1} \mb{X}^1  + \delta \mb{F}(\mb{X}^1, \Gamma^2), 
\\
\partial_t\mb{X}^2 &= \partial^2_{s^2} \mb{X}^2  + \delta \mb{F}(\mb{X}^2, \Gamma^1), 
\end{split}
\end{equation*}
which is subject to initial conditions $\mb{X}^1(\cdot, 0)$ and $\mb{X}^2(\cdot, 0)$ at the origin $t=0$. As a forcing term, we consider the Biot-Savart force $\mb{F}(\mb{X}^i, \Gamma^j)$ scaled by the factor $\delta=0.1$.

\begin{figure}
\begin{center}
\includegraphics[width=0.28\textwidth]{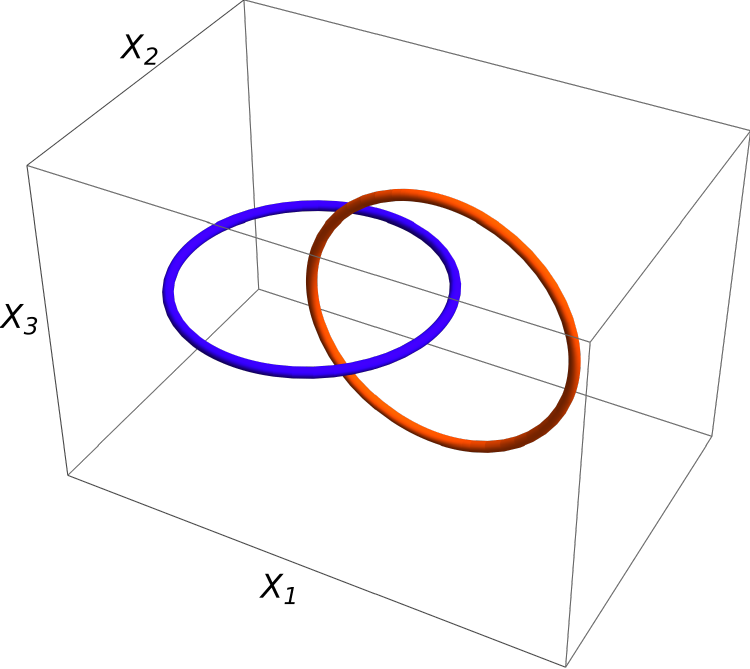}
\qquad \includegraphics[width=0.28\textwidth]{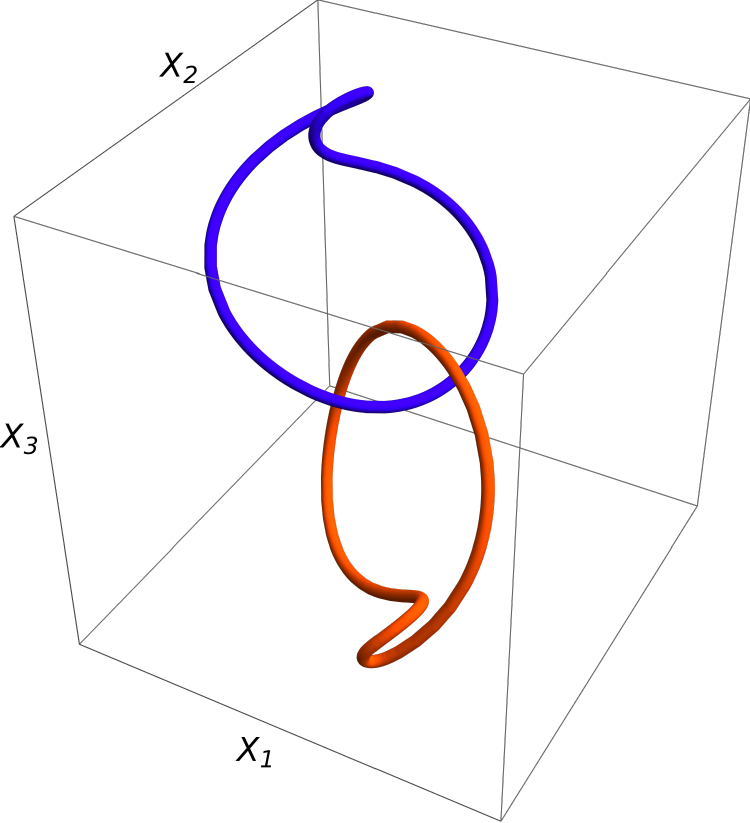}
\qquad \includegraphics[width=0.28\textwidth]{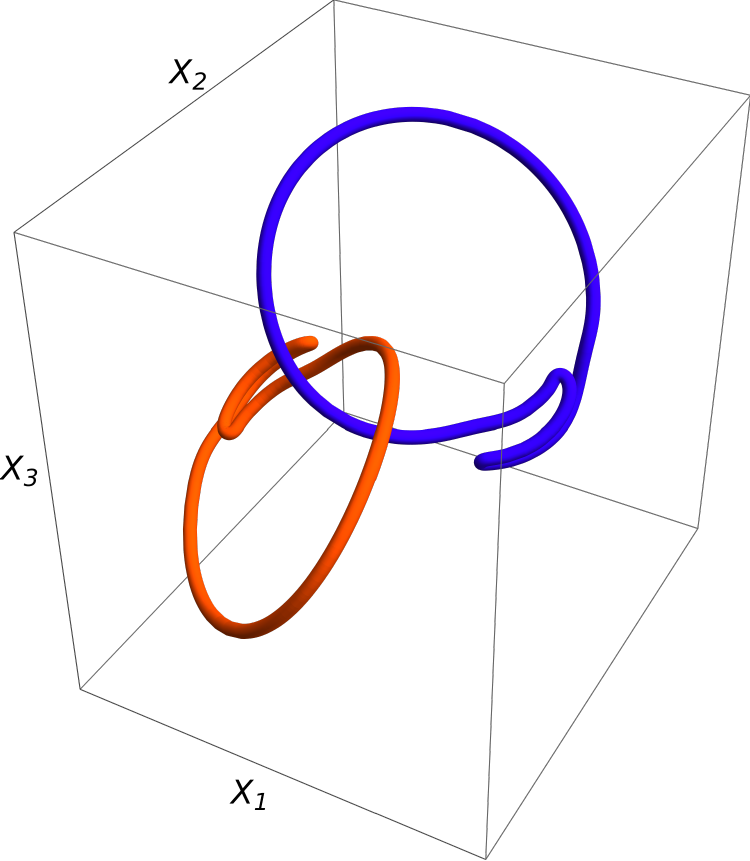}
\centerline{\scriptsize $t = 0$ \hskip 3truecm $t=0.031$ \hskip 3truecm $t=0.062$}

\smallskip

\includegraphics[width=0.28\textwidth]{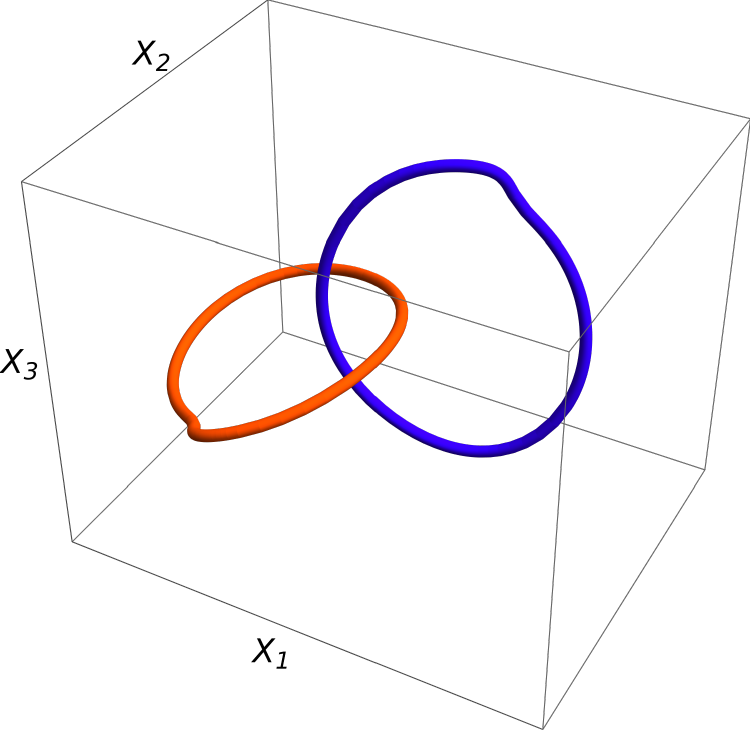}
\qquad \includegraphics[width=0.28\textwidth]{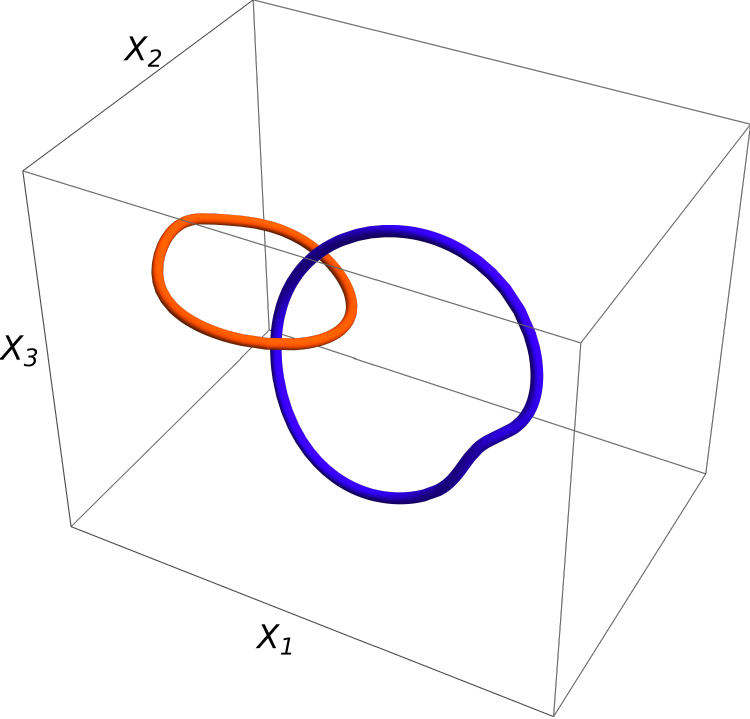}
\qquad \includegraphics[width=0.28\textwidth]{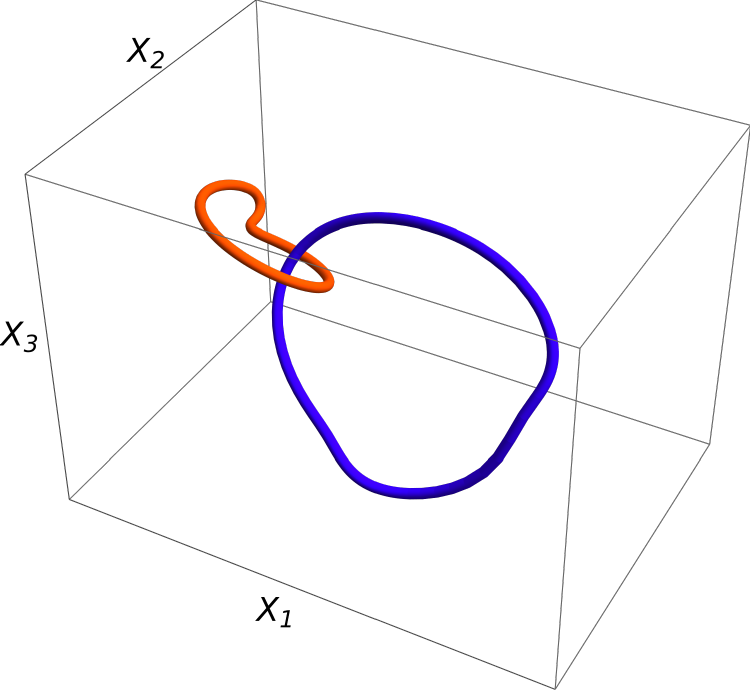}
\centerline{\scriptsize $t = 0.093$ \hskip 3truecm $t=0.124$ \hskip 3truecm $t=0.146$}

\caption{Evolution of the initial linked circles parameterized by (\ref{eq-crvs1-explicit}).}

\label{fig-knot-crvs1}
\end{center}
\end{figure}

\begin{figure}
\begin{center}
\includegraphics[width=0.28\textwidth]{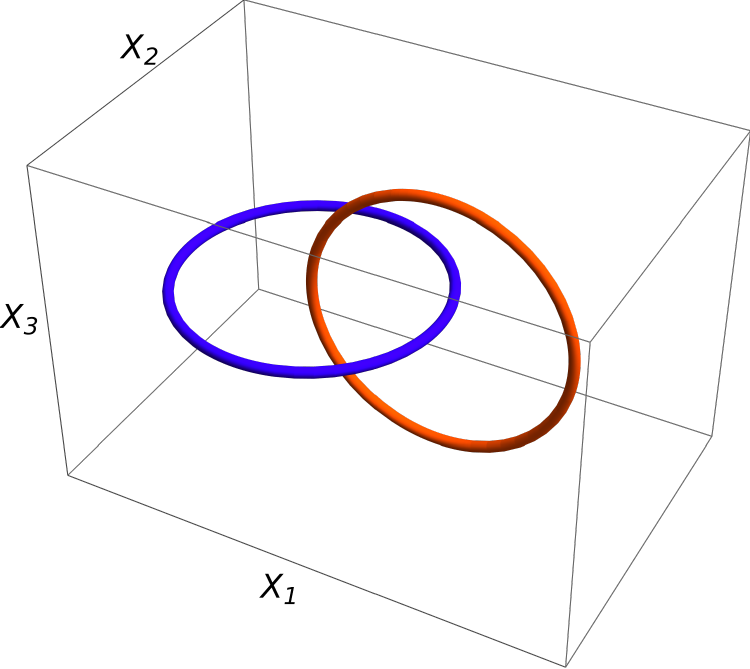}
\qquad \includegraphics[width=0.28\textwidth]{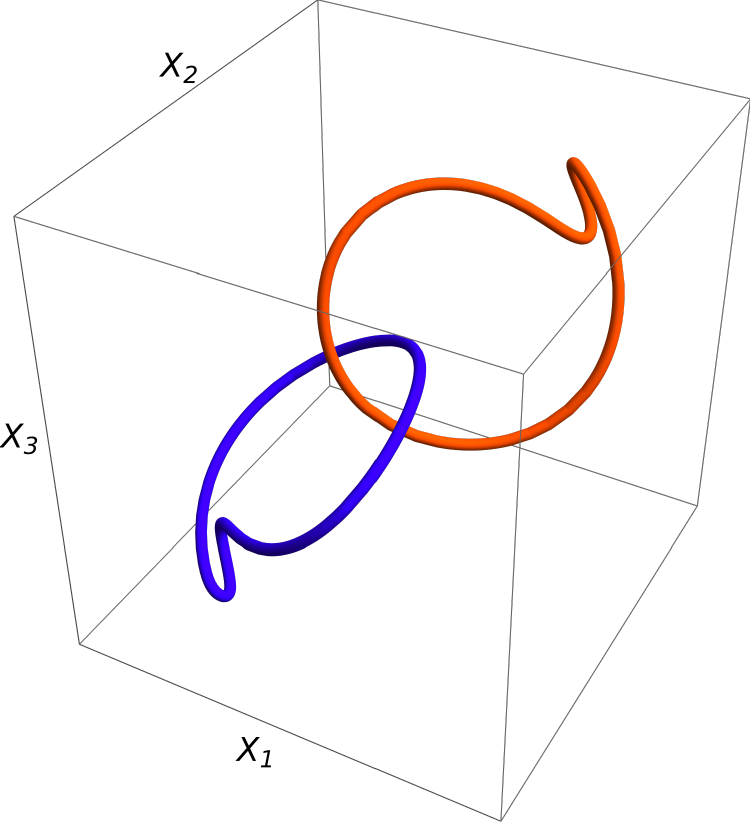}
\qquad \includegraphics[width=0.28\textwidth]{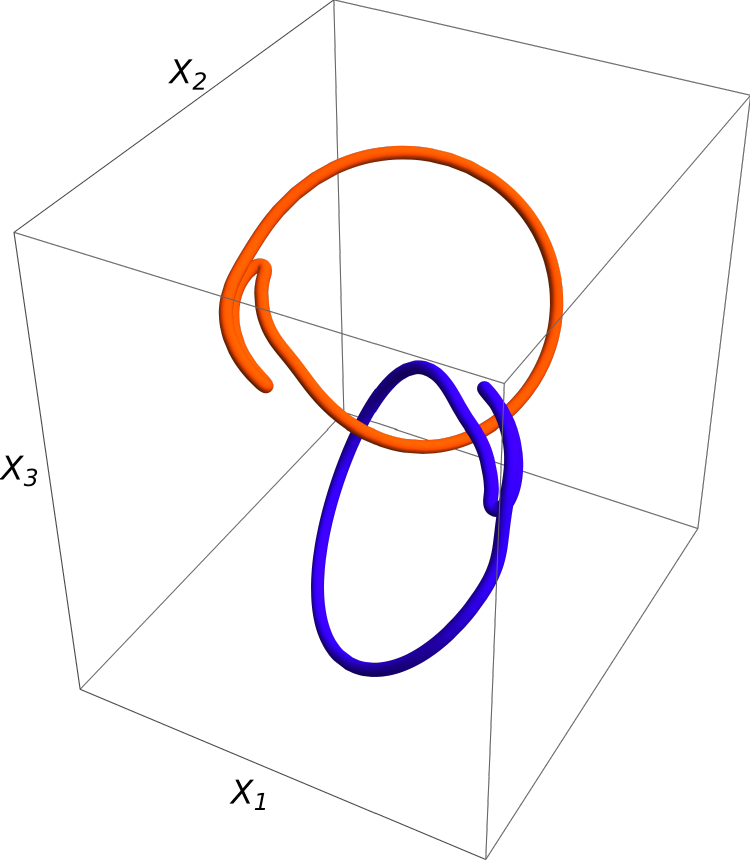}
\centerline{\scriptsize $t = 0$ \hskip 3truecm $t=0.031$ \hskip 3truecm $t=0.062$}

\smallskip

\includegraphics[width=0.28\textwidth]{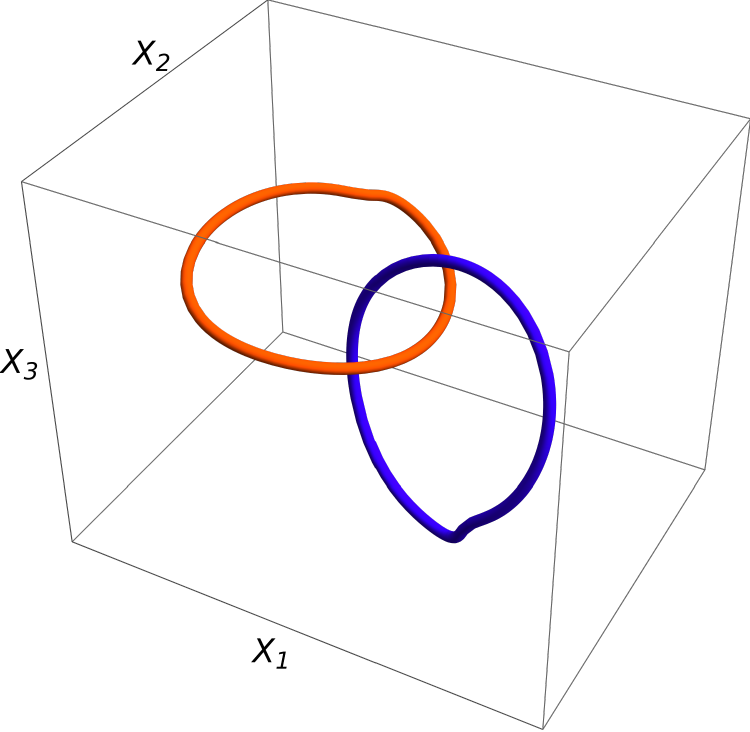}
\qquad \includegraphics[width=0.28\textwidth]{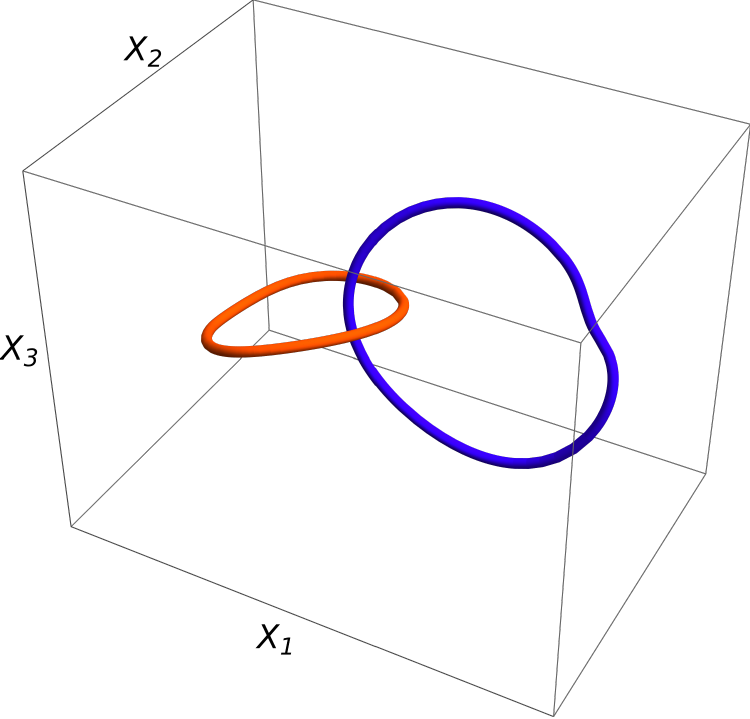}
\qquad \includegraphics[width=0.28\textwidth]{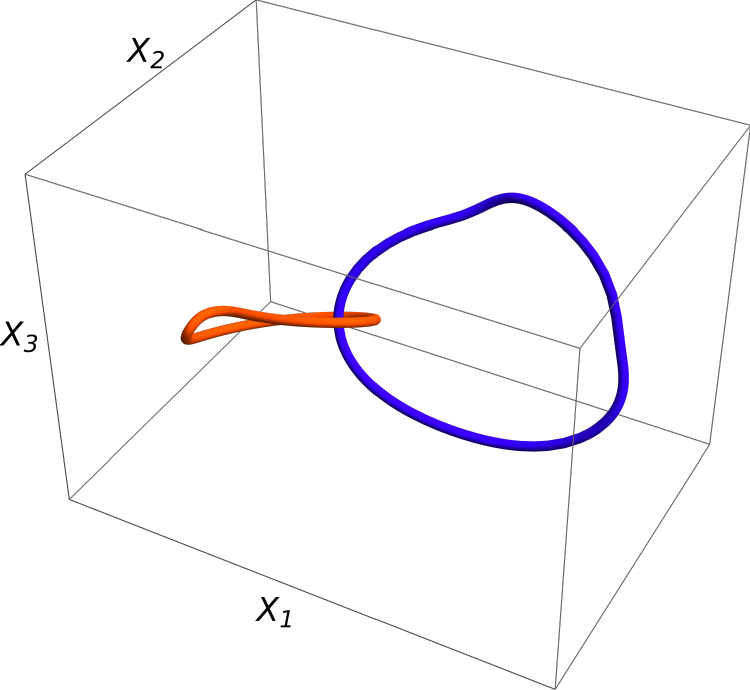}
\centerline{\scriptsize $t = 0.093$ \hskip 3truecm $t=0.124$ \hskip 3truecm $t=0.146$}

\caption{Evolution of the initial linked circles parameterized by (\ref{eq-crvs2-explicit}).}

\label{fig-knot-crvs2}
\end{center}
\end{figure}

\begin{figure}
\begin{center}
\includegraphics[width=0.42\textwidth]{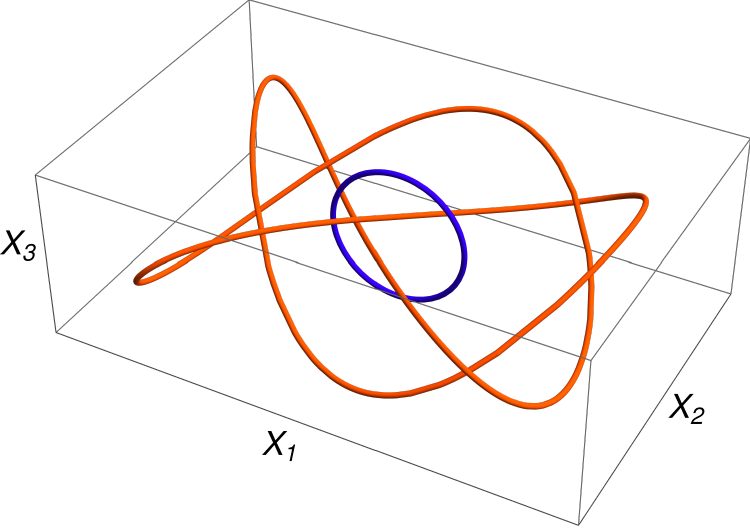}
\quad \includegraphics[width=0.42\textwidth]{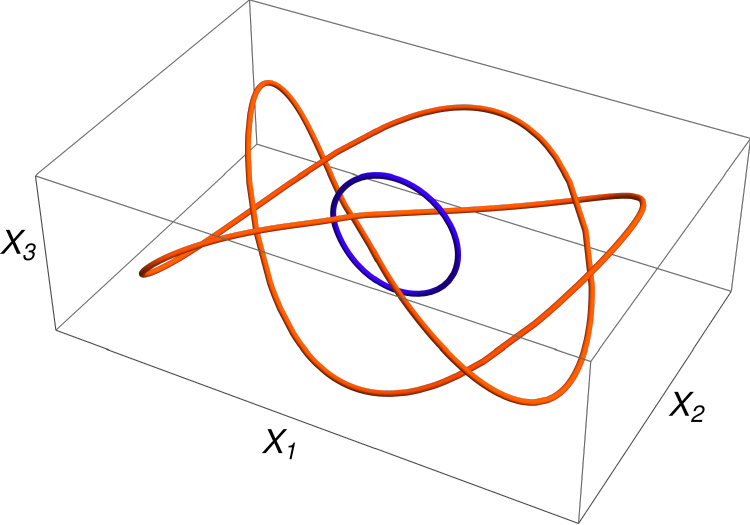}

\centerline{\scriptsize $t = 0$ \hskip 5truecm $t=0.077$}

\smallskip

\includegraphics[width=0.42\textwidth]{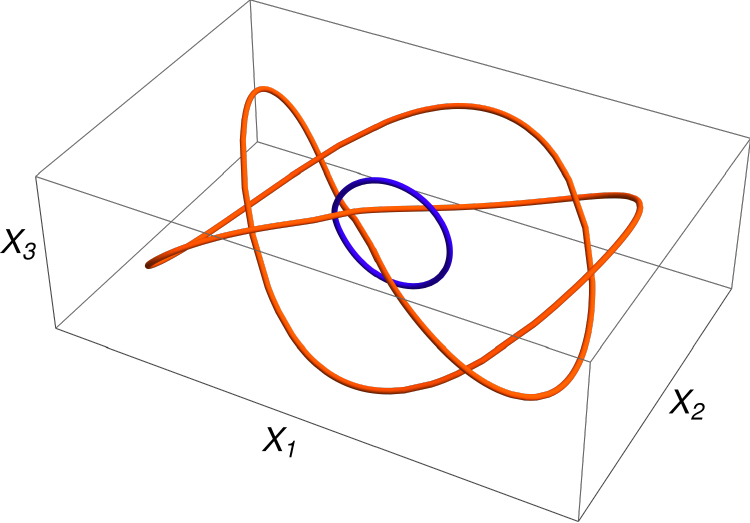}
\quad 
\includegraphics[width=0.42\textwidth]{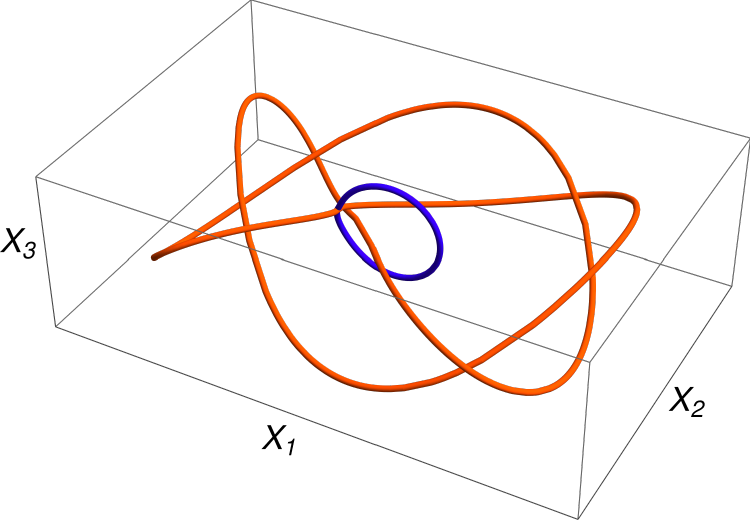}

\centerline{\scriptsize  $t=0.166$ \hskip 5truecm $t = 0.255$}

\smallskip

\includegraphics[width=0.42\textwidth]{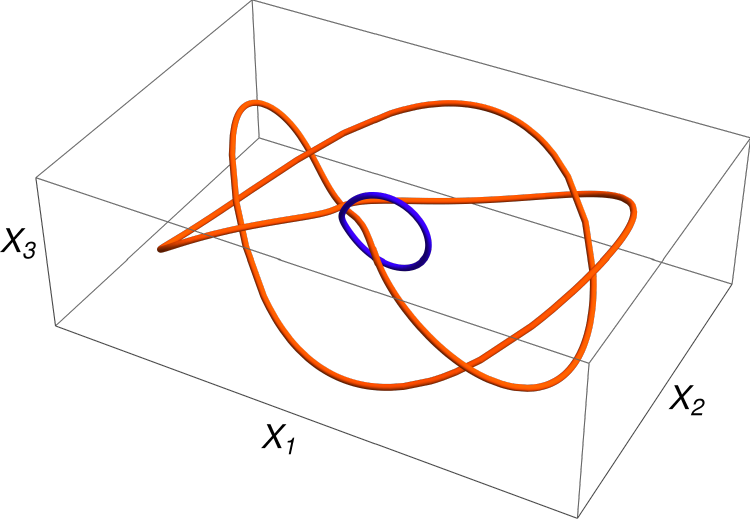}
\quad 
\includegraphics[width=0.42\textwidth]{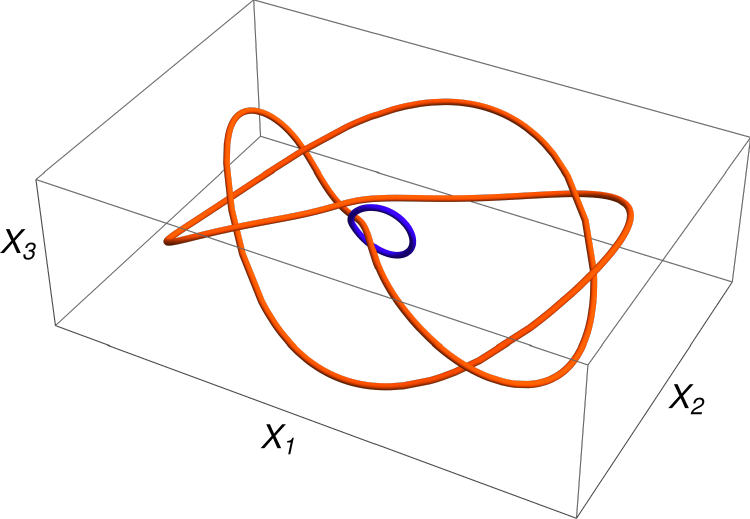}
\centerline{\scriptsize $t=0.344$ \hskip 5truecm $t=0.433$}

\caption{Evolution of the initial Listing's 8-curve linked with a circle (blue) (\ref{fig-knot-circles-ellipse}) a) and the external force given by the regularized Biot-Savart force given by (\ref{biot-savart-force}).}

\label{fig-knot-circle}
\end{center}
\end{figure}

\begin{figure}
\begin{center}
\includegraphics[width=0.42\textwidth]{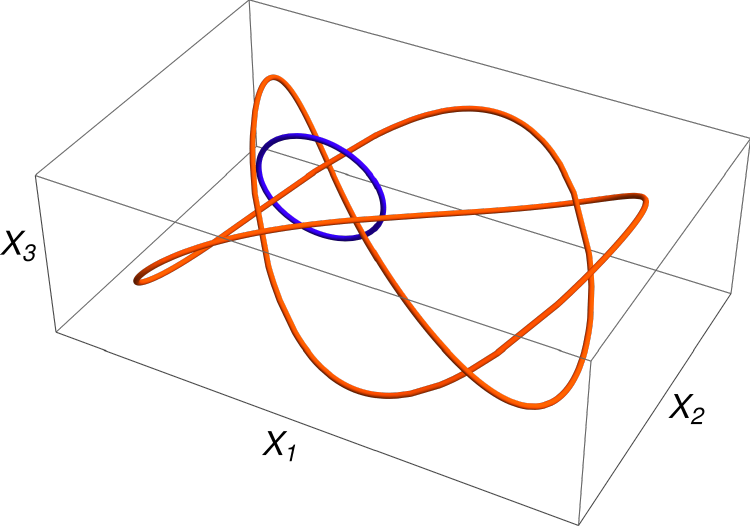}
\quad \includegraphics[width=0.42\textwidth]{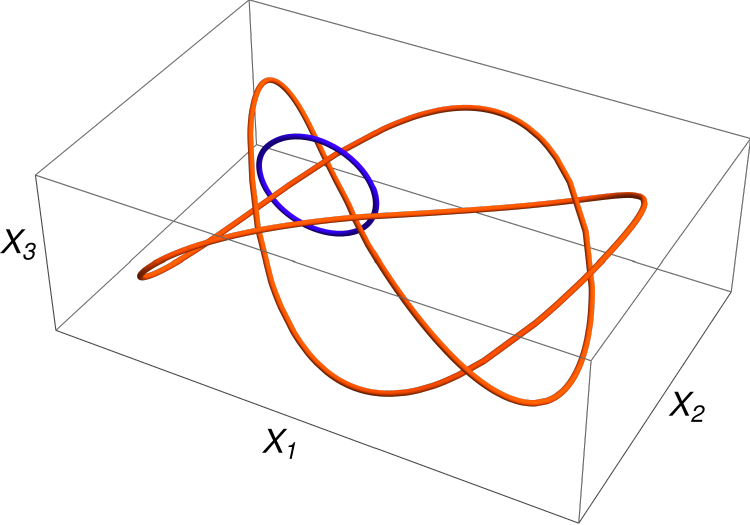}

\centerline{\scriptsize $t = 0$ \hskip 5truecm $t=0.049$}

\smallskip

\includegraphics[width=0.42\textwidth]{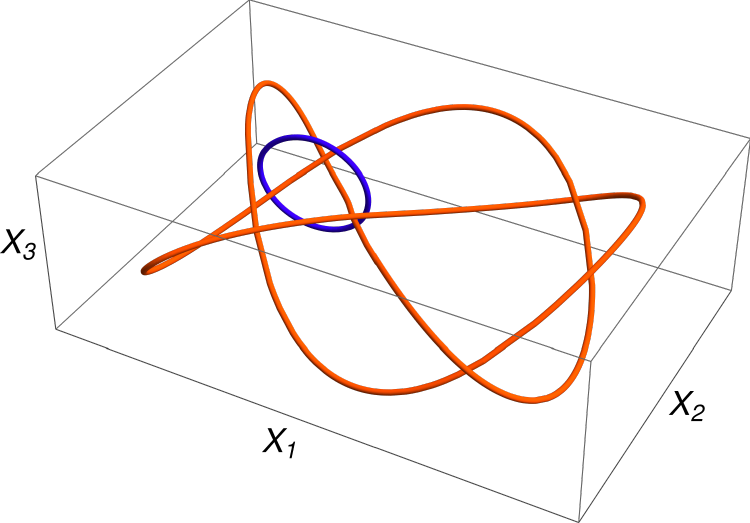}
\quad 
\includegraphics[width=0.42\textwidth]{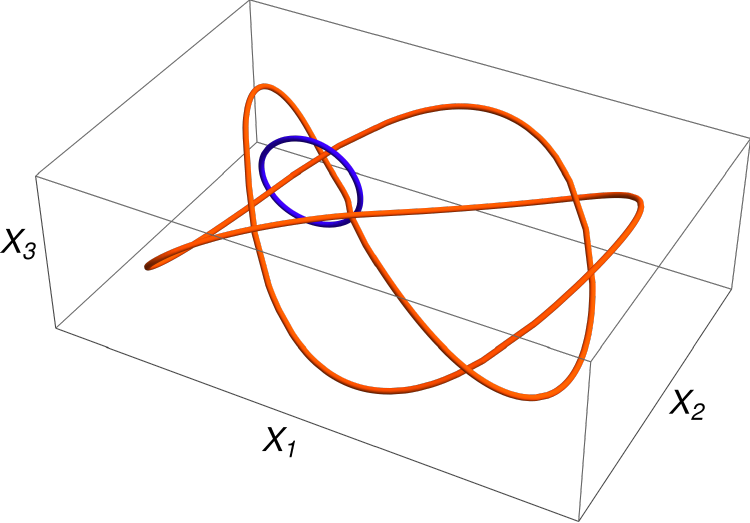}

\centerline{\scriptsize  $t=0.166$ \hskip 5truecm $t = 0.149$}

\smallskip

\includegraphics[width=0.42\textwidth]{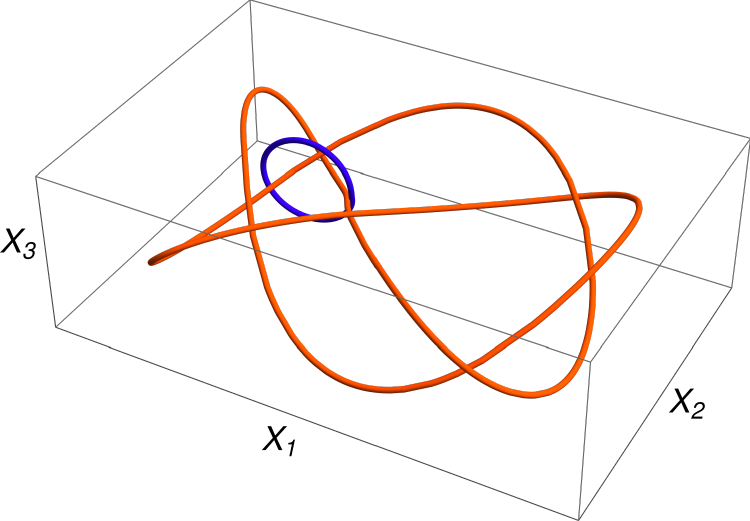}
\quad 
\includegraphics[width=0.42\textwidth]{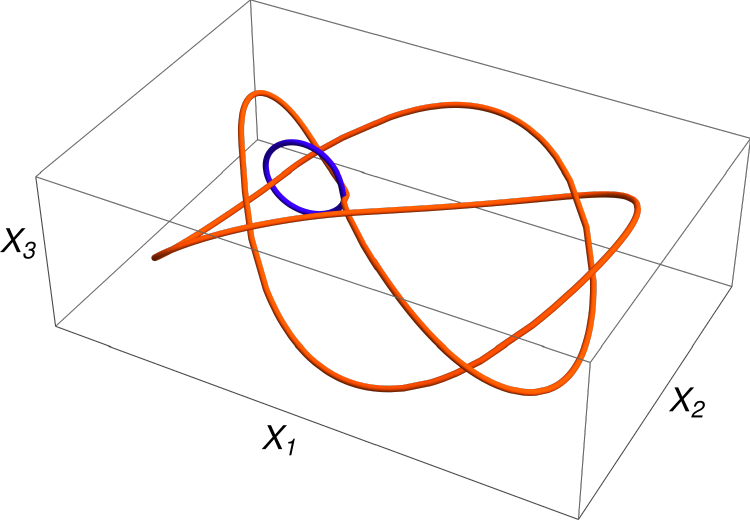}
\centerline{\scriptsize $t=0.198$ \hskip 5truecm $t=0.249$}

\caption{Evolution of the initial Listing's 8-curve linked with an ellipse (blue) (\ref{fig-knot-circles-ellipse}) b) and the external force given by the regularized Biot-Savart force given by (\ref{biot-savart-force}).}

\label{fig-knot-ellipse}
\end{center}
\end{figure}

In Fig.~\ref{fig-knot-crvs1} we present the time evolution of two linked circles $\Gamma^1$ and $\Gamma^2$ parameterized by (\ref{eq-crvs1-explicit}) with the linking number $link(\Gamma^1,\Gamma^2)=-1$. In Fig.~\ref{fig-knot-crvs2} we present the time evolution of two linked circles $\Gamma^1$ and $\Gamma^2$ parameterized by (\ref{eq-crvs2-explicit}) with the linking number $link(\Gamma^1,\Gamma^2)=1$. In Fig.~\ref{fig-knot-circle} we present the time evolution of the initial Listing's 8-curve $\Gamma^1$ linked with a circle $\Gamma^2$ parameterized by (\ref{eq-knot-circle1})  shown in Fig.~\ref{fig-knot-circles-ellipse} a) with the linking number $link(\Gamma^1,\Gamma^2)=0$. In Fig.~\ref{fig-knot-ellipse} we present the time evolution of the initial Listing's 8-curve $\Gamma^1$ linked with an ellipse circle $\Gamma^2$ shown in Fig.~\ref{fig-knot-circles-ellipse} b) with the linking number $link(\Gamma^1,\Gamma^2)=-2$.

\section{Conlusions}
In this article, we investigated a set of geometric evolution equations that describe the curvature-driven motion of a family of 3D curves along the normal and binormal directions. An evolving family of curves can interact in either local or non-local ways. In particular, we analyzed evolving pairs of closed linked curves that form knots in 3D. We utilized the direct Lagrangian method to solve the geometric flow of these interacting curves. We applied the abstract theory of nonlinear analytic semi-flows to prove the local existence, uniqueness, and continuation of classical H\"older smooth solutions for the system of nonlinear parabolic equations in question. Using the finite-volume method, we proposed an effective numerical method for solving the governing system of parabolic partial differential equations. Finally, we provided  multiple computational studies on the flow of linked curves.

\section*{Acknowledgments}
The first author received support from the project 21-09093S of the Czech Science Foundation.  The third author received support from the Slovak Research and Development Agency project APVV-20-0311.

\end{document}